*To my brother Mo*

# ANISOTROPIC QUADRATIC FORMS AND INVERSIVE GEOMETRY

NICHOLAS PHAT NGUYEN

Email: nicholas.pn@gmail.com

**Abstract**. We develop an inversive geometry for anisotropic quadradic spaces in analogy with the classical inversive geometry of a Euclidean plane. [1]

In this paper, we generalize the classical inversive geometry of a Euclidean plane to any vector space with an anisotropic quadratic form. Specifically, if $E$ is a vector space of finite dimension over a field $K$ of characteristic $\neq 2$, and we have an anisotropic symmetric bilinear form on $E$, then we can define certain transformations on the extended affine space $E \cup \{\infty\}$ (the space $E$ with an extra point $\infty$ called the point at infinity) that are analogous to the classical inversions of Euclidean plane geometry. Moreover, we can show that such inversive transformations have properties similar to the inversive transformations in the classical Euclidean plane.

An extension of classical inversive geometry to arbitrary anisotropic quadratic spaces depends on having the right notation and on the application of modern bilinear algebra. The extension applies only to anisotropic quadratic spaces and not to quadratic spaces in general, as the readers will find out why in this paper. Such a generalization has not been published and does not seem to be well-known, and therefore many readers may find the account in this paper useful and interesting, especially as the theory is rather elegant and could have new applications.

---

In particular, the readers can find some interesting applications of this paper in [1], [2], [3], and [4].

## Section 1. Anisotropic quadratic forms

Let $K$ be a field of characteristic $\neq 2$. In this paper, we will consider anisotropic quadratic forms defined for vector spaces of finite dimension over $K$.

Anisotropic quadratic forms are those quadratic forms $h$ where $h(X) = 0$ if and only if $X = 0$. We will denote the symmetric bilinear form (scalar product) associated with $h$ by using the dot product notation. Specifically, we will write $h(\mathrm{X}) = X.X$, where $(X, Y) \mapsto X.Y$ is the associated symmetric bilinear form. We have $2X.Y = h(X + Y) - h(X) - h(Y)$. This symmetric bilinear form is obviously regular or non-degenerate because $h$ is anisotropic.

Anisotropic quadratic forms occur naturally in many mathematical contexts.

- Multiplication on $K$ naturally gives us an anisotropic form of dimension 1.
- If $E$ is a quadratic field extension of $K$, then the norm map from $E$ to $K$ is naturally an anisotropic quadratic form of dimension 2. When $K$ is a finite field, this is essentially the only anisotropic quadratic space of dimension > 1.
- If $E$ is a quaternion algebra defined over $K$, then the reduced norm is a quadratic form of dimension 4 that is anisotropic if and only if $E$ is a division algebra. (If the reduced norm is isotropic, then $E$ is isomorphic to a matrix algebra.)
- Given an anisotropic quadratic form, the extension by scalars of that form to a field $L$ over $K$ is also anisotropic when $L$ is a finite extension of odd degree over $K$, or when $L$ is the rational function field $K(T)$. See [5] at chapter 2, section 5.
- Let $q$ be quadratic form over a number field $K$. If $q$ is anisotropic over one of the completions of $K$, then clearly $q$ must be anisotropic over $K$. In fact, the converse is also true. It was a landmark result in number theory and the beginning of the local-global principle in modern mathematics. Specifically, the Hasse-Minkowski theorem says that if $q$ is isotropic over every completion of a number field $K$, then it is also isotropic over $K$.
- If $K$ is an ordered field, the standard scalar product

$$\boldsymbol{x}.\boldsymbol{y} = x_1 y_1 + \ldots + x_n y_n \text{ where } \boldsymbol{x} = (x_1, \ldots, x_n) \text{ and } \boldsymbol{y} = (y_1, \ldots, y_n)$$

is always anisotropic.

The fact that the standard scalar product $\boldsymbol{x}.\boldsymbol{y} = x_1 y_1 + \ldots + x_n y_n$ is anisotropic for any $K^n$ is a characteristic property of ordered field, as discovered by Artin and Schreier (a field can be ordered if and only if it is formally real).

In addition to the standard scalar product, any positive or negative definite quadratic form over an ordered field would of course be anisotropic. As an application of the ideas outlined in this paper, the interested reader may want to consider the case of a positive definite quadratic form in a real Euclidean space. The theory outlined in this paper would provide a coherent and extensive theory of inversive geometry in that context, something that people are still trying to work out with some isolated results.

Note that if $K$ is an algebraically closed field, then there is no anisotropic quadratic space of dimension $\geq 2$. If $K$ is a $p$-adic number field (an extension of finite degree of $\mathbf{Q}_p$), then there is no anisotropic quadratic space of dimension $\geq 5$.

If $K$ is a finite field, then there is no anisotropic quadratic space of dimension $\geq 3$, because any homogeneous polynomial of degree 2 in 3 or more variables will automatically have a non-trivial zero over $K$, by the Chevalley-Warning theorem. A function field of an algebraic curve over an algebraically closed field also has the same property by Tsen's theorem. So there is also no anisotropic quadratic space of dimension $\geq 3$ in the case of a function field of an algebraic curve over an algebraically closed field.

## Section 2. The projective space of cycles for an anisotropic quadratic form

For an anisotropic quadratic space $E$ over $K$, consider the set $F$ of all functions $p$ from $E$ to $K$ of the form $p(\boldsymbol{X}) = a\boldsymbol{X}.\boldsymbol{X} + \boldsymbol{b}.\boldsymbol{X} + c$, where $\boldsymbol{X}$ and $\boldsymbol{b}$ are vectors in $E$, and $a$ and $c$ are elements in the field $K$.[2]

The set $F$ of all such functions can naturally be endowed with the structure of a $K$-vector space of dimension $\dim E + 2$, being parametrized by the vector $\boldsymbol{b}$ and the elements $a$ and $c$. We will refer to a function $p$ where the coefficient $a$ is non-zero as a 2-cycle, and to a function $p$ where the coefficient $a = 0$ but $\boldsymbol{b} \neq \boldsymbol{0}$ as a linear cycle or 1-cycle. We refer to the constant functions $p(\boldsymbol{X}) = c$ as 0-cycles.

Aside from the natural structure of a vector space over $K$, we can also endow $F$ with a symmetric bilinear form $<_\,,_>$ as follows.

Given $p = a\boldsymbol{X}.\boldsymbol{X} + \boldsymbol{b}.\boldsymbol{X} + c$ and $p^* = a^*\boldsymbol{X}.\boldsymbol{X} + \boldsymbol{b^*}.\boldsymbol{X} + c^*$, define

$$<p,p^*> \;=\; \boldsymbol{b}.\boldsymbol{b^*} - 2ac^* - 2a^*c.$$

This scalar product is clearly symmetric and bilinear. Moreover, it is non-degenerate, because $F$ is isometric to the sum of $E$ and an Artinian plane (also known as a hyperbolic plane). We will refer to this fundamental scalar product on $F$ as the cycle pairing or cycle product.

For our purpose, we consider only functions $p \neq 0$, and moreover we regard all non-zero scalar multiples of the same function as equivalent because they have the same zero set in $E$ (if any). In other words, let $Z(p) = \{\boldsymbol{x} \in E \text{ such that } p(\boldsymbol{x}) = 0\}$. Then $Z(p) = Z(up)$ where $u$ is any non-zero element of $K$.

---

[2] If $E$ has finite dimension $n$, and we regard the coefficients of the vector $\boldsymbol{X}$ in a chosen basis of $E$ as variables, the function $p$ is at most a second-degree polynomial in n variables. Because the field $K$ has 3 or more elements, such a function $p$ is a zero function if and only if the parameters $a$, $\boldsymbol{b}$, and $c$ are all zero.

Accordingly, we are led naturally to the projective space $\hat{F}$ consisting of all equivalent classes of such non-zero functions, or what is essentially the same, the set of all 1-dimensional subspaces of $F$. As a matter of convenient notation, for each letter denoting an element or subset of $F$, we will put a hat over that letter to denote the corresponding element or subset of $\hat{F}$ under the natural correspondence that sends elements of $F$ to their projective equivalence classes.

The projective space $\hat{F}$ can be decomposed into the following disjoint subsets:

(i) A subset of all 2-cycles $\hat{p}$ where $p = a\boldsymbol{X}.\boldsymbol{X} + \boldsymbol{b}.\boldsymbol{X} + c$ with coefficient $a \neq 0$.

By completing the square, we can rewrite each 2-cycle as $p = a\boldsymbol{X}.\boldsymbol{X} + \boldsymbol{b}.\boldsymbol{X} + + c = a(\boldsymbol{X} + \boldsymbol{b}/2a).\ (\boldsymbol{X} + \boldsymbol{b}/2a) - as$, where $s = (\boldsymbol{b}.\boldsymbol{b} - 4ac)/4a^2 = <p, p>/4a^2$. By geometric analogy, we will refer to the point $(-\mathbf{b}/2\text{a})$ in $E$ as the center of the 2-cycle $p$, and the number $s$ as the size of the 2-cycle. In case $K = \mathbf{R}$, the number $s$ is equal to the square of the radius $r$ if $\text{Z}(p)$ is an ordinary geometric circle or sphere.

Note that all 2-cycles have centers, although some 2-cycles may not have any zeros in $E$. The 2-cycles with size $s = 0$ are of course isotropic 2-cycles, i.e., 2-cycles $p$ such that $<p, p> = 0$.

(ii) A projective subspace corresponding to all 1-cycles and 0-cycles $p = \boldsymbol{b}.\boldsymbol{X} + c$. The constant functions $p(\boldsymbol{X}) = c$ represent an exceptional point of this subset, which we will denote by the infinity symbol $\infty$.

There is a natural bijection between isotropic 2-cycles in the projective space $\hat{F}$ and the points of $E$. For any point in $E$, take the set of all 2-cycles with zero radius centered at that point. These 2-cycles are all scalar multiples of each other, and map to a unique projective element in $\hat{F}$.

Using standard language for quadratic spaces, we say that $u$ is orthogonal to $v$ if $<u, v> = 0$. Since orthogonal relations are obviously unchanged when we replace each vector by a scalar multiple, such relations are consistent with projective equivalence and it makes sense to talk about orthogonal elements in the projective space $\hat{F}$. For any subset $S$ of $F$ (or $\hat{F}$) we denote by $S^{\#}$ the subspace of elements of $F$ (or $\hat{F}$) that are orthogonal under the cycle pairing to all elements in $S$. A vector $u$ of $F$ is said to be isotropic if $<u, u> = 0$ and non-isotropic or anisotropic otherwise. This concept is clearly compatible with the projective equivalence, so we can similarly refer to isotropic and non-isotropic points in the projective space $\hat{F}$.

The proposition below summarizes the basic orthogonal relationships in $F$.

**<u>Proposition 2.1</u>**:

*(a)* *Any non-zero constant function (0-cycle) is orthogonal to itself and all 1-cycles, but is not orthogonal to any 2-cycle.*

*(b)* *Two 1-cycles $p = \boldsymbol{b}.\boldsymbol{X} + c$ and $q = \boldsymbol{d}.\boldsymbol{X} + e$ have scalar product $<p, q> = \boldsymbol{b}.\boldsymbol{d}$. Consequently, $p$ and $q$ are orthogonal in $F$ if and only if $\boldsymbol{b}$ and $\boldsymbol{d}$ are orthogonal in $E$. A 1-cycle is therefore non-isotropic.*

*(c)* *A 2-cycle $p$ with center $C$ is orthogonal to a 1-cycle $L$ if and only if $L(C) = 0$, i.e. the affine hyperplane defined by the equation $L(\boldsymbol{x}) = 0$ passes through $C$.*

*(d)* *For 2-cycles $p = a(\boldsymbol{X} + \boldsymbol{b}/2a).\ (\boldsymbol{X} + \boldsymbol{b}/2a) - as$ and $q = d(\boldsymbol{X} + \boldsymbol{e}/2d).\ (\boldsymbol{X} + \boldsymbol{e}/2d) - dt$, the scalar product of $p$ and $q$ are given the following formula:*

$$<p, q> = 2ad(s + t - w), \text{ where } w = (\boldsymbol{b}/2a - \boldsymbol{e}/2d).\ (\boldsymbol{b}/2a - \boldsymbol{e}/2d)$$

*(e)* *The only isotropic cycles in $F$ are the 0-cycles and the 2-cycles of size zero, i.e. functions $p(\boldsymbol{X})$ of the form $a(\boldsymbol{X} + \boldsymbol{b}/2a).\ (\boldsymbol{X} + \boldsymbol{b}/2a)$.*

***Proof***. Statements (a) and (b) are clear.

For statement (c), let $p = a\boldsymbol{X}.\boldsymbol{X} + \boldsymbol{b}.\boldsymbol{X} + c = a(\boldsymbol{X} + \boldsymbol{b}/2a).\ (\boldsymbol{X} + \boldsymbol{b}/2a) - as$ be a 2-cycle with center $C = -\boldsymbol{b}/2a$. Let $q = \boldsymbol{e}.\boldsymbol{X} + f$ be a 1-cycle. Then $<p, q> = \boldsymbol{b}.\boldsymbol{e} - 2af = -2a.q(C) = 0$ if and only if $q(C) = 0$.

For statement (d), let:

$$p = a\boldsymbol{X}.\boldsymbol{X} + \boldsymbol{b}.\boldsymbol{X} + c = a(\boldsymbol{X} + \boldsymbol{b}/2a).\ (\boldsymbol{X} + \boldsymbol{b}/2a) - as, \text{ where } s = (\boldsymbol{b}.\boldsymbol{b} - 4ac)/4a^2$$

$$q = d\boldsymbol{X}.\boldsymbol{X} + \boldsymbol{e}.\boldsymbol{X} + f = d(\boldsymbol{X} + \boldsymbol{e}/2d).\ (\boldsymbol{X} + \boldsymbol{e}/2d) - dt, \text{ where } t = (\boldsymbol{e}.\boldsymbol{e} - 4df)/4d^2 ,$$

then by definition $<p, q> = \boldsymbol{b}.\boldsymbol{e} - 2af - 2dc$.

Because $w = (\boldsymbol{b}/2a - \boldsymbol{e}/2d).\ (\boldsymbol{b}/2a - \boldsymbol{e}/2d) = \boldsymbol{b}.\boldsymbol{b}/4a^2 + \boldsymbol{e}.\boldsymbol{e}/4d^2 - 2\boldsymbol{b}.\boldsymbol{e}/4ad$, we have, after cancelling out the terms $\boldsymbol{b}.\boldsymbol{b}/4a^2$ and $\boldsymbol{e}.\boldsymbol{e}/4d^2$:

$$s + t - w = - c/a - f/d + 2\boldsymbol{b}.\boldsymbol{e}/4ad, \text{ or}$$

$$2ad(s + t - w) = \boldsymbol{b}.\boldsymbol{e} - 2af - 2dc = <p, q>$$

Statement (e) follows from (a) through (d). For a 2-cycle $p$, $<p, p> = \boldsymbol{b}.\boldsymbol{b} - 4ac = 2a^2(2s - 0) = 0$ if and only if $s = 0$. ■

## Section 3. Zero sets

Let $V$ be the set of isotropic points in $\hat{F}$. We just saw that $V$ consists of: (i) the exceptional point $\infty$ representing the 0-cycles, and (ii) points represented by 2-cycles of size zero, which can be naturally identified with the points of $E$. Accordingly, we can identify $V$ with $E \cup \{\infty\}$. With this identification, we can think of $V$ as an extension of the affine space $E$, and we will refer to the exceptional point $\infty$ of $V$ as the point at infinity of the affine space $E$.

We also have the following description for the zero set Z($p$) of any cycle $p$ in $F$. Recall that we define $Z(p)$ as the set $\{\boldsymbol{x} \in E \text{ such that } p(\boldsymbol{x}) = 0\}$.

**Proposition 3.1**: *Let p be a non-zero cycle of F.*

*(a) The orthogonal complement $(p)^{\#}$ of $(p)$ is a linear subspace of co-dimension 1 (i.e., a hyperplane). If p is non-isotropic, then the hyperplane $(p)^{\#}$ is regular under the scalar product induced by the cycle pairing, and F is isometric to the direct sum of $(p)$ and $(p)^{\#}$.*

*(b) If p is a 1-cycle, then the intersection $\widehat{(p)}^{\#} \cap V$ contains the following: (i) the exceptional point $\infty$ and (ii) all those points of E that lie on the affine hyperplane $Z(p)$.*

*Because of this fact, we will, as a matter of convenience, consider the exceptional point $\infty$ as a zero of any 1-cycle p, i.e., that Z(p) includes the point $\infty$.*

*(c) If p is a 2-cycle, then the intersection $\widehat{(p)}^{\#} \cap V$ contains precisely all those points of E that lie on the zero set Z(p) of p. Z(p) can be empty, in which case $\widehat{(p)}^{\#}$ does not intersect V.*

***Proof*.** Statement (a) follows from the basic facts about regular symmetric bilinear spaces. Statement (b) is an immediate consequence of Proposition 2.1 (a) and (c). Statement (b) implies that the exceptional point $\infty$ belongs to the orthogonal complement $\widehat{(p)}^{\#}$ of any 1-cycle $p$.

In most cases, the essential set that we look at is not really the zero set Z($p$), but the orthogonal complement $\widehat{(p)}^{\#}$ of $p$ and particularly its intersection with $V$. For a 1-cycle $p$, such a set includes not just the points on the affine hyperplane Z($p$), but also the exceptional point $\infty$ of $V$.

For statement (c), let $p$ is a 2-cycle. The 0-cycles are not orthogonal to $p$, and so the zero set Z($p$) of $p$ does not contain the point $\infty$. An isotropic 2-cycle $q$ is orthogonal to $p$ if and only if (using the same notation as in Proposition 2.1), the expression $s + t - w = s - w = 0$. But we have $aw = a(\mathbf{b}/2a - \mathbf{e}/2d).(\mathbf{b}/2a - \mathbf{e}/2d) = p(-\mathbf{e}/2d) + as$. Therefore $s = w$ if and only if $as = aw$ or $p(-\mathbf{e}/2d) = 0$. This means the center of the isotropic 2-cycle $q$ belongs to the zero set Z($p$) of $p$. ■

The zero set Z($p$) can of course be empty, e.g., consider the case $s < 0$ when $K = \mathbf{R}$. In that case, the orthogonal complement $\widehat{(p)}^{\#}$ does not intersect with $V$. We can say more about the zero set Z($p$) for 2-cycles $p$.

**<u>Proposition 3.2:</u>** *Let p be a 2-cycle in F. For convenience, we identify the zero set Z(p) in E with the intersection $V \cap \widehat{(p)}^{\#}$ by virtue of Proposition 3.1.*

*(a) Z(p) consists of exactly one point if and only if p is an isotropic 2-cycle.*

(b) *If p is a non-isotropic 2-cycle, then Z(p) is non-empty if and only if* $<p, p>$ *is represented by the anisotropic quadratic form on E. In that case, Z(p) contains as many points as an affine hyperplane in the extended plane* $E \cup \{\infty\}$.

***Proof.*** (a) Proposition 3.1(c) implies that if $p$ is an isotropic 2-cycle, then Z($p$) consists of exactly one point. Moreover, if $p$ is a non-isotropic 2-cycle, then the orthogonal complement $(p)^{\#}$ is regular, and therefore either has no isotropic vectors, or contains an Artinian plane (with two linearly independent isotropic vectors). Accordingly, Z($p$) consists of exactly one point if and only if the 2-cycle $p$ is isotropic.

(b) Let $p$ be a non-isotropic 2-cycle, $p(\boldsymbol{X}) = a\boldsymbol{X}.\boldsymbol{X} + \boldsymbol{b}.\boldsymbol{X} + c = a(\boldsymbol{X} + \boldsymbol{b}/2a).\,(\boldsymbol{X} + \boldsymbol{b}/2a) - as$,

where $s = (\boldsymbol{b}.\boldsymbol{b} - 4ac)/4a^2 = <p, p>/4a^2$ is non-zero.

Z($p$) is non-empty if and only if there is a vector $\boldsymbol{z}$ in E satisfying the equation

$$(\boldsymbol{z} + \boldsymbol{b}/2a).\,(\boldsymbol{z} + \boldsymbol{b}/2a) = s.$$

That means $s$ can be represented by the anisotropic quadratic form on $E$. Because $<p, p> = 4a^2s$, the number $<p, p>$ can also be represented by the same quadratic form on $E$ exactly when $s$ is.

Now let $L$ be the 1-cycle $(\boldsymbol{z} + \boldsymbol{b}/2a).\boldsymbol{X}$. We have $<L, L> = (\boldsymbol{z} + \boldsymbol{b}/2a).\,(\boldsymbol{z} + \boldsymbol{b}/2a) = s$. The subspace ($p$) is isometric to ($L$) because $<p, p> = 4a^2s$ and so the numbers $<p, p>$ and $<L, L>$ differ by a square factor. Given the isometry between the subspaces ($p$) and ($L$), Witt's extension theorem applied to the quadratic space $F$ says that there is an orthogonal transformation of $F$ mapping $p$ to $L$.[3] The same transformation would map $(p)^{\#}$ to $(L)^{\#}$, and the intersection $(p)^{\#} \cap$ V to the intersection $(L)^{\#} \cap V$, which consists of the affine hyperplane Z($L$) and $\infty$. ■

It is worthwhile to note the following geometrical interpretation of the cycle pairing defined on the vector space $F$ in the classical case of a Euclidean space $\mathbf{R}^n$.

---

[3] For a statement and proof of Witt's extension theorem (also known as Witt's cancellation theorem), see, e.g., [1], theorem 3.9 or [5], chapter 7.

**Proposition 3.3**: *For $E = \boldsymbol{R}^n$, let p and q be 2-cycles with sizes s and t respectively.*

(a) *Assume that both s and $t \geq 0$. In other words, Z(p) and Z(q) are ordinary spheres in $\boldsymbol{R}^n$. Then $<p, q> = 0$ if and only if the associated spheres are orthogonal in the classical sense.*

(b) *If both s and t are $< 0$, then $<p, q> \neq 0$.*

(c) *Assume that $s \geq 0$ but $t < 0$. Let $q^*$ be the 2-cycle obtained by replacing t by the number $-t$ in the expression for q. We have $<p, q> = 0$ if and only if the intersection of Z(p) and Z($q^*$) lies in an affine hyperplane passing through the center of Z($q^*$).*

***Proof***. Note that $<p, q> = 0$ is equivalent to $s + t - w = 0$ in light of Proposition 2.1(d). When $s$ and $t$ are both $> 0$,

$s$ = square of the radius of Z($p$)

$t$ = square of the radius of Z($q$)

$w$ = square of the distance between the centers of Z($p$) and Z($q$).

The equation $s + t - w = 0$ is then a simple restatement of the classical meaning of orthogonal spheres in the Euclidean space $\mathbf{R}^n$.

If both $s$ and $t$ are $< 0$, then $s + t - w < 0$ (as $w$ is by definition always $> 0$, being the square of the distance between two center points). Accordingly, $<p, q>$ is never zero in this case.

If $s \geq 0$ but $t < 0$, then rewrite $s + t - w = 0$ as $s = -t + w$. That means the square of the radius of Z($p$) is equal to the square of the radius of Z($q^*$) plus the square of the distance between the center points of Z($p$) and Z($q^*$). ($q$ and $q^*$ have the same center points.) That happens if and only if the hyperplane that passes through the intersection of Z($p$) and Z($q^*$) also passes through the center point of Z($q^*$). When $n = 2$, that is the same as saying that Z($p$) cuts Z($q^*$) in a diameter of Z($q^*$). ■

## <u>Section 4. Group action on isotropic points</u>.

We consider the group $O(F)$ of orthogonal transformations of the vector space $F$ that leaves invariant the cycle pairing defined on $F$. In other words, $O(F)$ consists of all linear transformations $L \in GL(F)$ such that $<L(p), L(q)> \; = \; <p, q>$ for any cycles p, q in $F$.

We denote by $\mathcal{G}$ the corresponding group of projective orthogonal transformations on F, i.e., it is the subgroup of $PGL(F)$ induced by orthogonal transformations in $O$(F).

The most basic orthogonal transformations in $O(F)$ are the hyperplane reflections defined by non-isotropic cycles. Given a non-isotropic cycle $p$ in $F$, we can express $F$ as the direct sum of $(p)$ and its orthogonal complement $(p)^{\#}$. Any cycle $\boldsymbol{x}$ of F is then equal to a unique sum $\boldsymbol{x} = \boldsymbol{u} + \boldsymbol{v}$, where $\boldsymbol{u} \in (p)$ and $\boldsymbol{v} \in (p)^{\#}$. The reflection $\Re$ across the orthogonal hyperplane $(p)^{\#}$ is the linear transformation: $\boldsymbol{x} = \boldsymbol{u} + \boldsymbol{v} \;\rightarrow - \boldsymbol{u} + \boldsymbol{v}$.

Expressed in terms of $p$, the orthogonal reflection $\Re$ has the following well-known formula:

$$\Re(x) = x - 2\,\frac{<x,p>}{<p,p>}\,p$$

Note that $< p, p > \neq 0$ because $p$ is a non-isotropic cycle. As the reader can readily verify from the above formula, $\Re$ sends $p$ to $-p$, and for any cycle $\boldsymbol{v}$ orthogonal to $p$, $\Re(\boldsymbol{v}) = \boldsymbol{v}$.

The theorem of Cartan – Dieudonne says that the group $O(F)$ is generated by hyperplane reflections. More precisely, each orthogonal transformation in $O(F)$ can be expressed as product of at most $\dim(F) = \dim(E) + 2$ hyperplane reflections defined by non-isotropic vectors of $F$.

According to Proposition 2.1, the set $V$ of isotropic cycles of $\hat{F}$ can be naturally identified with $E \cup \{\infty\}$. Each orthogonal transformation in $\mathcal{G}$ maps $V$ to itself, as the image of an isotropic vector is clearly another isotropic vector. We therefore have a natural group action of $\mathcal{G}$ on $V$.

**Theorem 4.1**: *The natural group action of $\mathcal{G}$ on V is faithful and doubly transitive. Moreover, if a regular hyperplane $\widehat{H}$ intersects V in a non-empty set W, then a transformation in $\mathcal{G}$ that leaves all points of W fixed must either be the identity transformation or the orthogonal reflection across the hyperplane $\widehat{H}$.*

***Proof***. To say that the action of $\mathcal{G}$ on $V$ is doubly transitive means that we can always find an orthogonal projective transformation that moves any two distinct points **u**, **w** of $V$ to any other two distinct points **u'**, **w'** of V. Indeed, let $M$ be the 2-dimensional subspace of $F$ spanned by {**u**, **w**} and $N$ the subspace spanned by {**u'**, **w'**}. Both $M$ and $N$ are isometric to an Artinian (hyperbolic) plane. By Witt's extension theorem, there is an orthogonal projective transformation in $\mathcal{G}$ that maps **u** to **u'** and **w** to **w'**.

To say that the action of $\mathcal{G}$ on $V$ is faithful means that the only transformation in $\mathcal{G}$ leaving all elements of $V$ invariant is the identity transformation. Recall that we can identify $V$ with $E \cup \{\infty\}$, where $\infty$ is identified with the 0-cycles, and each point of $E$ is identified with the isotropic 2-cycles centered at that point. Let $\boldsymbol{w}$ be the zero vector in $E$, and let $\boldsymbol{u_1}, \boldsymbol{u_2}, \ldots, \boldsymbol{u_n}$ be an orthogonal basis in $E$, where $n = \dim(E)$. Then under the identification of $E \cup \{\infty\}$ and $V$, the points $\boldsymbol{w}, \boldsymbol{u_1}, \boldsymbol{u_2}, \ldots, \boldsymbol{u_n}$, and $\infty$ correspond to the following linear basis of the vector space $F$:

$$v_0 = \boldsymbol{X.X}$$

$$v_1 = \boldsymbol{X.X} - 2\boldsymbol{u_1.X} + \boldsymbol{u_1.u_1}$$

…………..

$$v_n = \boldsymbol{X.X} - 2\boldsymbol{u_n.X} + \boldsymbol{u_n.u_n}$$

$$v_{n+1} = 1$$

In addition to the points $\boldsymbol{w}, \boldsymbol{u_1}, \boldsymbol{u_2}, \ldots, \boldsymbol{u_n}$, and $\infty$, we can also find another point $\boldsymbol{z}$ in $V$ whose expression as linear combination of $v_0, v_1, \ldots, v_n$, and $v_{n+1}$ has no zero coefficient. Indeed, let $U = \boldsymbol{u_1.u_1} + \boldsymbol{u_2.u_2} + \ldots + \boldsymbol{u_n.u_n}$. We can always choose an integer $j$ such that both numbers $(-j - n + 2)$ and $(j^2 - 1)\,\boldsymbol{u_1.u_1} + U$ are nonzero. Now consider $\boldsymbol{z} = (j\boldsymbol{u_1} + \boldsymbol{u_2} + \ldots + \boldsymbol{u_n})$ in $E$. As an element of $V$, it is the point corresponding to the 2-cycle $\boldsymbol{X.X} - 2(j\boldsymbol{u_1} + \boldsymbol{u_2} + \ldots + \boldsymbol{u_n}).\boldsymbol{X} + (j^2 - 1) + U = (-j - n + 2)v_0 + jv_1 + v_2 + \ldots + v_n + ((j^2 - 1)\, \boldsymbol{u_1.u_1} + U)\, v_{n+1}$ in $F$.

From linear algebra, we know that the $(n + 3)$ points of $V$ corresponding to $\boldsymbol{w}, \boldsymbol{u_1}, \boldsymbol{u_2}, \ldots, \boldsymbol{u_n}, \infty$ and $z$ form a projective frame of $\hat{F}$. This implies that any projective transformation of $\hat{F}$ which leaves these points fixed must be the identity transformation itself. Hence, the action of $\mathcal{G}$ on V is certainly faithful.

Now suppose a regular hyperplane $\hat{H}$ intersects $V$ in a non-empty set $W$. According Proposition 3.2, if a regular hyperplane $\hat{H}$ in $\hat{F}$ intersects $V$ in a non-empty set $W$, then there is a projective orthogonal transformation of $\hat{F}$ transforming $\hat{H}$ into $(\hat{L})^{\#}$ where $L = \boldsymbol{b.X}$ is a linear cycle with $\boldsymbol{b.b}$ non-zero. We can pick $\boldsymbol{b}$ to be part of an orthogonal basis of E, and so without loss of generality we can assume that $\hat{H}$ is the projective hyperplane orthogonal to $\mathbf{u_1}$.

In that case, $H$ can be regarded as a space of cycles in $(n - 1)$ variables and the set $W$ includes $\mathbf{w}, \mathbf{u_2}, \ldots, \mathbf{u_n}, \infty$. Any orthogonal transformation leaving all of the points of $W$ fixed must also leave all vectors in $H$ invariant. The only orthogonal transformations with that property are the identity transformation and the orthogonal reflection across $H$. ■

If $p$ is a 1-cycle, then the orthogonal reflection defined by $p$ clearly maps 1-cycless to 1-cycles and 2-cycles to 2-cycles. If $p$ is a non-zero 2-cycle, the formula $\Re(x) = x - 2\dfrac{\langle x,p\rangle}{\langle p,p\rangle}\, p$ shows that the orthogonal reflection defined by $p$ maps a 1-cycle $x$ to another 1-cycle if and only if $\langle x, p\rangle = 0$, i.e., if the zero set of the 1-cycle $x$ passes through the center of $p$.

<u>**Proposition 4.2**</u>. *If p is a non-isotropic 2-cycle, the orthogonal reflection* $\Re(x) = x - 2\dfrac{\langle x,p\rangle}{\langle p,p\rangle}\, p$ *maps a non-isotropic 2-cycle x to a 1-cycle if and only if Z(x) passes through the center of p.*

***Proof***. Because the formula for $\Re(x)$ is linear in $x$ and independent of any linear factor of $p$, we can assume that both $x$ and $p$ have leading coefficient 1. In that case, the orthogonal reflection $\Re$ maps a 2-cycle $x$ to a 1-cycle if and only if $2\langle x, p\rangle = \langle p, p\rangle$.

Let $p = \boldsymbol{X.X} - 2\boldsymbol{b.X} + c$ (with center $\mathbf{b}$), and $x = \boldsymbol{X.X} + \boldsymbol{d.X} + e$. We have $<p, p> = 4\boldsymbol{b.b} - 4c$, and $<x, p> = -2\boldsymbol{d.b} - 2e - 2c$.

Accordingly, $<p, p> - 2<x, p> = 4\boldsymbol{b.b} - 4c + 4\boldsymbol{d.b} + 4e + 4c = 4x(\boldsymbol{b})$. That means we have $2<x, p> = <p, p>$ if and only if $x(\boldsymbol{b}) = 0$, i.e., $\boldsymbol{b}$ is a zero of the 2-cycle $x$. ■

## Section 5. Reflections and inversions – Inversive transformations

Recall that we start with an anisotropic quadratic space $E$. On the extended affine space $E \cup \{\infty\}$ obtained by adjoining an extra point $\infty$ called the point at infinity, we can define certain transformations that are analogous to the classical inversions of Euclidean plane geometry.

The zero set of any 1-cycle $L(\boldsymbol{x}) = \boldsymbol{b.x} + c$ where $\boldsymbol{b} \neq \mathbf{0}$ is an affine hyperplane in $E$. The corresponding affine reflection is the mapping that sends $\infty$ to itself, and any finite point $\boldsymbol{x}$ to a finite point $\boldsymbol{x}' = \boldsymbol{x} - 2\boldsymbol{b}\dfrac{L(\boldsymbol{x})}{\boldsymbol{b.b}}$. Note that $\boldsymbol{b.b}$ is a non-zero number because the space $E$ is anisotropic and $\boldsymbol{b} \neq \mathbf{0}$ by hypothesis.

Likewise, for any non-isotropic 2-cycle $a(\boldsymbol{x} - \boldsymbol{b}).(\boldsymbol{x} - \boldsymbol{b}) - as$, where $a$ and $s$ are non-zero elements of $F$, we can define an "inversion" mapping the space $E \cup \{\infty\}$ to itself as follows:

- $\boldsymbol{b} \leftrightarrow \infty$,
- for any finite point $\boldsymbol{x} \neq \boldsymbol{b}$, $\boldsymbol{x} \rightarrow$ the finite point $\boldsymbol{x}'$ collinear with $\boldsymbol{b}$ and $\boldsymbol{x}$, and such that $(\boldsymbol{x} - \boldsymbol{b}).(\boldsymbol{x}' - \boldsymbol{b}) = s$.

The collinear condition and dot product equation are equivalent to $\boldsymbol{x}' - \boldsymbol{b} = (\boldsymbol{x} - \boldsymbol{b})\dfrac{s}{(\boldsymbol{x}-\boldsymbol{b}).(\boldsymbol{x}-\boldsymbol{b})}$. Note that $(\boldsymbol{x} - \boldsymbol{b}).(\boldsymbol{x} - \boldsymbol{b})$ is a non-zero number because $(\boldsymbol{x} - \boldsymbol{b}) \neq \mathbf{0}$ and the quadratic space $E$ is anisotropic by hypothesis. This shows that inversion is a well-defined transformation regardless of whether or not the 2-cycle $a(\boldsymbol{x} - \boldsymbol{b}).(\boldsymbol{x} - \boldsymbol{b}) - as$ has any zero in the space $E$.

Any non-zero scalar multiple of a non-isotropic cycle defines the same mapping, and any reflection or inversion so defined is a transformation of $E \cup \{\infty\}$ that is its own inverse. If **x** and **x'** are mapped to each other by such a transformation, we will call these points inverse points or conjugate points relative to the non-isotropic cycle that defines the transformation.

We denote by *Inv(E)* the group of transformations of $E \cup \{\infty\}$ generated by such reflections and inversions. We will refer to the elements of *Inv(E)* as inversive transformations.

Let $p$ be a non-isotropic cycle in $F$ and let $\mathcal{R}$ be the hyperplane reflection defined by $p$. We will use the same letter $\mathcal{R}$ to denote the above orthogonal transformation in $O(F)$ and the projective transformation in $\mathcal{G}$ induced by the same reflection. Because the meaning will be clear from the context, there should be little risk of confusion.

The following theorem describes the action of those reflections on $V$.

**<u>Theorem 5.1</u>**:

(a) *If $p$ is a 1-cycle $\boldsymbol{b}.\boldsymbol{X} + c$, then the action of $\mathcal{R}$ on $V$ is the same as the affine reflection defined by $p$.*

(b) *If $p$ is a non-isotropic 2-cycle $(\boldsymbol{X} - \boldsymbol{b}).(\boldsymbol{X} - \boldsymbol{b}) - s$ with $s \neq 0$, then the action of $\mathcal{R}$ on $V$ is the same as the inversion with respect to the 2-cycle $p$.*

**Proof**. We will prove (a) and (b) by showing that in each case the action of $\mathcal{R}$ on $V$ coincides with such an affine reflection or cycle inversion, which we denote by the letter $L$.

Recall that the reflection $\mathcal{R}$ is the orthogonal transformation with the formula

$$\mathcal{R}(q) = q - 2\frac{\langle q,p \rangle}{\langle p,p \rangle}\, p$$

We begin with statement (a). When $p$ is the 1-cycle $\boldsymbol{b}.\boldsymbol{X} + c$, $\langle p, p\rangle = \boldsymbol{b}.\boldsymbol{b}$. We know that the exceptional point $\infty$ is orthogonal to $p$. That means $\mathcal{R}(\infty) = \infty = L(\infty)$.

For any finite point $\boldsymbol{w}$ in $V$ represented by an isotropic 2-cycle $q = \boldsymbol{X}.\boldsymbol{X} - 2\boldsymbol{w}.\boldsymbol{X} + \boldsymbol{w}.\boldsymbol{w}$, we note first that $L(\boldsymbol{w})$ by definition is the point $\boldsymbol{w}' = \boldsymbol{w} - \boldsymbol{b}\,\dfrac{2p(\boldsymbol{w})}{\boldsymbol{b}.\boldsymbol{b}}$.

We have $<q, p> = -2\boldsymbol{b}.\boldsymbol{w} - 2c = -2p(\boldsymbol{w})$. Hence,

$$\mathcal{R}(q) \;=\; q - 2\,\frac{<q,p>}{<p,p>}\,p \;=\; q\; - 2\,\frac{-2p(\boldsymbol{w})}{\boldsymbol{b}.\boldsymbol{b}}\,(\boldsymbol{b}.\boldsymbol{X} + c)$$

$$= \boldsymbol{X}.\boldsymbol{X} - 2(\boldsymbol{w} - \boldsymbol{b}\,\frac{2p(\boldsymbol{w})}{\boldsymbol{b}.\boldsymbol{b}}).\boldsymbol{X} + \boldsymbol{w}.\boldsymbol{w} + 4c\,\frac{p(\boldsymbol{w})}{\boldsymbol{b}.\boldsymbol{b}}$$

$$= \boldsymbol{X}.\boldsymbol{X} - 2\,\boldsymbol{w}'.\boldsymbol{X} + \boldsymbol{w}.\boldsymbol{w} + 4c\,\frac{p(\boldsymbol{w})}{\boldsymbol{b}.\boldsymbol{b}}$$

Now observe that:

$$\boldsymbol{w}'.\boldsymbol{w}' = \boldsymbol{w}.\boldsymbol{w} + \frac{4p(\boldsymbol{w})p(\boldsymbol{w})}{\boldsymbol{b}.\boldsymbol{b}} - \frac{4p(\boldsymbol{w})\boldsymbol{b}.\boldsymbol{w}}{\boldsymbol{b}.\boldsymbol{b}} = \boldsymbol{w}.\boldsymbol{w} + 4\,\frac{p(\boldsymbol{w})}{\boldsymbol{b}.\boldsymbol{b}}\,(p(\boldsymbol{w}) - \boldsymbol{b}.\boldsymbol{w}) = \boldsymbol{w}.\boldsymbol{w} + 4\,\frac{p(\boldsymbol{w})}{\boldsymbol{b}.\boldsymbol{b}}c$$

Substituting the above in the expression for $\mathcal{R}(q)$, we have $\mathcal{R}(q) = \boldsymbol{X}.\boldsymbol{X} - 2\boldsymbol{w}'.\boldsymbol{X} + \boldsymbol{w}'.\boldsymbol{w}'$.

What the above expression means is that the finite point $\boldsymbol{w}$ in $V$ is mapped to $\boldsymbol{w}'$ under the reflection $\mathcal{R}$. But we already see that $\boldsymbol{w}' = \mathrm{L}(\boldsymbol{w})$ is the image of $\boldsymbol{w}$ under the reflection across the affine hyperplane $\mathrm{Z}(p)$. Accordingly, $\mathcal{R}$ has the same action on $V$ as the affine reflection $L$. This proves statement (a).

We now consider statement (b). Let $p$ be the non-isotropic 2-cycle $(\boldsymbol{X} - \boldsymbol{b}).\,(\boldsymbol{X} - \boldsymbol{b}) - s = \boldsymbol{X}.\boldsymbol{X} - 2\boldsymbol{b}.\boldsymbol{X} + + \boldsymbol{b}.\boldsymbol{b} - s$ with center $\boldsymbol{b}$ and size $s \neq 0$.

The inversion $L$ with respect to the 2-cycle $p$ is by definition the mapping:

- $\boldsymbol{b} \leftrightarrow \infty$,
- for any finite point $\boldsymbol{x} \neq \boldsymbol{b}$, $\boldsymbol{x} \rightarrow$ the finite point $\boldsymbol{x}'$ collinear with $\boldsymbol{b}$ and $\boldsymbol{x}$, and that satisfies the dot product equation $(\boldsymbol{x} - \boldsymbol{b}).(\boldsymbol{x}' - \boldsymbol{b}) = s$.

For any vector $q$, we have $\Re(q) = q - 2\,\frac{<q,p>}{<p,p>}\,p$. From the defining formula for the cycle pairing, it is easy to verify that:

$$<p, p> = (-2\boldsymbol{b}).(-2\boldsymbol{b}) - 4(\boldsymbol{b}.\boldsymbol{b} - s) = 4\boldsymbol{b}.\boldsymbol{b} - 4\boldsymbol{b}.\boldsymbol{b} + 4s = 4s.$$

For the constant function $q(\boldsymbol{X}) = s$, $<q, p> = -2s$.

Hence, if $q$ is the constant function $q(\boldsymbol{X}) = s$,

$$\Re(q) = q - 2\,\frac{<q,p>}{<p,p>}\,p = s + p = (\boldsymbol{X} - \boldsymbol{b}).(\boldsymbol{X} - \boldsymbol{b})$$

But the zero 2-cycle $(\boldsymbol{X} - \boldsymbol{b}).(\boldsymbol{X} - \boldsymbol{b}) = \boldsymbol{X}.\boldsymbol{X} - 2\boldsymbol{b}.\boldsymbol{X} + \boldsymbol{b}.\boldsymbol{b}$ is exactly the image of $\boldsymbol{b}$ under the identification $V \leftrightarrow E \cup \{\infty\}$. Hence $\boldsymbol{b} \leftrightarrow \infty$ under the action of both $\Re$ and $L$.

For any finite point $\boldsymbol{m} \neq \boldsymbol{b}$, take the zero 2-cycle $q(\boldsymbol{X}) = \boldsymbol{X}.\boldsymbol{X} - 2\boldsymbol{m}.\boldsymbol{X} + \boldsymbol{m}.\boldsymbol{m}$ corresponding to $\boldsymbol{m}$ under the identification $V \leftrightarrow E \cup \{\infty\}$. We know that the image $\Re(q)$ of $q$ will be projectively equivalent to another isotropic 2-cycle $q'(\boldsymbol{X}) = \boldsymbol{X}.\boldsymbol{X} - 2\boldsymbol{m}'.\boldsymbol{X} + \boldsymbol{m}'.\boldsymbol{m}'$ associated with some point $\boldsymbol{m}'$. Our task is to determine $\boldsymbol{m}'$.

In this case, $<q, p> = 4\boldsymbol{m}.\boldsymbol{b} - 2\boldsymbol{m}.\boldsymbol{m} - 2\boldsymbol{b}.\boldsymbol{b} + 2s = 2(s - (\boldsymbol{m} - \boldsymbol{b}).(\boldsymbol{m} - \boldsymbol{b}))$. Accordingly,

$$\Re(q) = q - 2\,\frac{<q,p>}{<p,p>}\,p = q - \left(1 - \frac{(\boldsymbol{m}-\boldsymbol{b}).(\boldsymbol{m}-\boldsymbol{b})}{s}\right)p$$

If $q' = \boldsymbol{X}.\boldsymbol{X} - 2\boldsymbol{m}'.\boldsymbol{X} + \boldsymbol{m}'.\boldsymbol{m}'$ is equivalent to $\Re(q)$, then we must have

$$\lambda q' = q - \left(1 - \frac{(\boldsymbol{m}-\boldsymbol{b}).(\boldsymbol{m}-\boldsymbol{b})}{s}\right)p \quad \text{for some non-zero scalar } \lambda.$$

Because we are working with $\boldsymbol{m} \neq \boldsymbol{b}$, the scalar product $(\boldsymbol{m} - \boldsymbol{b}).(\boldsymbol{m} - \boldsymbol{b})$ is a non-zero number. To simplify our expression, let $\alpha$ be the number in $K$ such that

$$1/\alpha = \frac{(\boldsymbol{m}-\boldsymbol{b}).(\boldsymbol{m}-\boldsymbol{b})}{s} \quad \text{or equivalently} \quad \alpha(\boldsymbol{m} - \boldsymbol{b}).(\boldsymbol{m} - \boldsymbol{b}) = s.$$

We then have $\lambda q' = q - (1 - 1/\alpha)p$.

Comparing the respective coefficients of $\boldsymbol{X.X}$ and $\boldsymbol{X}$ from each side, we have the equations:

$$\lambda = 1/\alpha$$

$$\lambda \boldsymbol{m'} = \boldsymbol{m} - (1 - 1/\alpha)\boldsymbol{b}$$

Because $\alpha\lambda = 1$, we have:

$$\boldsymbol{m'} = \alpha\boldsymbol{m} - (\alpha - 1)\boldsymbol{b}$$

or $$(\boldsymbol{m'} - \boldsymbol{b}) = \alpha(\boldsymbol{m} - \boldsymbol{b})$$

But the last equation implies that (i) $\boldsymbol{m'}$ is collinear with $\boldsymbol{m}$ and $\boldsymbol{b}$, and (ii) $(\boldsymbol{m'} - \boldsymbol{b}).(\boldsymbol{m} - \boldsymbol{b}) = \alpha(\boldsymbol{m} - \boldsymbol{b}).(\boldsymbol{m} - \boldsymbol{b}) = s$. These are exactly the conditions for $\boldsymbol{m'}$ and $\boldsymbol{m}$ to be conjugate or inverse points with respect to the 2-cycle $p$.

To sum up, the orthogonal reflection $\mathfrak{R}$ defined by the non-isotropic 2-cycle $p$ has the same action on $V$ as the inversion $L$ defined by $p$. This proves statement (b), and our proof of Theorem 5.1 is complete. ■

Because the 1-cycles and non-isotropic 2-cycles $p$ in $F$ for which a reflection or inversion on the space of $E \cup \{\infty\}$ can be defined are exactly the non-isotropic vectors for which we can define reflections in $F$, the above theorem shows that all inversive transformations in *Inv(E)* come from orthogonal transformations in $\mathcal{G}$. Because orthogonal transformations map hyperplanes to hyperplanes while leaving $V$ stable, inversive transformations map a cycle zero set to another cycle zero set. Moreover, such action is orthogonal in the sense that it maps orthogonal zero sets to orthogonal zero sets.

The classical study of pencils and bundles of circles can best be viewed as the study of 1-dimentional and 2-dimensional subspaces of $\hat{F}$. For example, consider the classical case of the Euclidean plane $\mathbf{R}^2$. The corresponding vector space $F$ of cycles has linear dimension $2 + 2 = 4$ and the projective space $\hat{F}$ has dimension $n + 1 = 3$. Let $\hat{H}$ be a 2-dimensional subspace of $\hat{\mathrm{F}}$.

Then $H$ is a hyperplane in $F$, and the orthogonal complement of $H$ has dimension 1. That means we can always describe $\widehat{H}$ as the set of cycles orthogonal to some $\hat{p}$ in $\widehat{F}$. The traditional classification of such bundles defined by $H$ into parabolic, hyperbolic and elliptic types follows naturally from the above general consideration based on the character of the cycle $p$ (isotropic, non-isotropic with non-empty zero set, non-isotropic with empty zero set).

The traditional classification of pencils into 3 different types also has a natural explanation. Let $\hat{L}$ be a 1-dimensional projective subspace in $\hat{F}$. Looking at the scalar product on the 2-dimensional vector space $L$ induced by the cycle pairing, we have 3 possibilities:

(i) $L$ is singular. In this case, $L$ cannot be totally isotropic because the isotropic variety $V$ contains no subspace of dimension 1 or more. Hence, its radical is 1-dimensional and the all cycles in the pencil $\hat{L}$ are orthogonal to the cycle representing that radical. Consequently, all cycles in the pencil $\hat{L}$ have a common zero in $V$.

(ii) $L$ is regular and anisotropic (the pencil $\hat{L}$ does not intersect $V$). In this case, its orthogonal complement $L^{\#}$ is also regular and must contain isotropic cycles. ($F$ is the direct sum of $L$ and $L^{\#}$ and $F$ contains isotropic cycles.) All cycles in $L$ are orthogonal to these isotropic cycles, which means that the points of $V$ represented by these isotropic cycles are the common zeros of all the cycles in $L$.

(iii) $L$ is regular and isotropic. In this case, $L$ is an Artinian plane. So $L$ has two linearly independent and non-orthogonal isotropic cycles, and the pencil $\hat{L}$ intersects $V$ at exactly two distinct points.[4] As we will see in further discussion below, the points in $V$ corresponding to these isotropic cycles are conjugate or inverse points relative to any cycle in the pencil $L$.

When $E$ is finite-dimensional, Witt's extension theorem implies that if $L$ has two independent isotropic vectors (i.e., $L$ is an Artinian plane), then its orthogonal complement $L^{\#}$ must have none (anisotropic), and vice versa.

---

[4] Over a field of characteristic $\neq 2$, a plane with a regular symmetric bilinear form is isotropic if and only if it is an Artinian plane (also known as a hyperbolic plane).

## Section 6. Conjugate points – Stereographic projections

Recall that a pair of points in $E \cup \{\infty\}$ are said to be conjugate or inverse with respect to a non-isotropic cycle $p$ if they are mapped to each other under the reflection or inversion defined by $p$. Theorem 5.1 shows that two points are conjugate with respect to $p$ if and only if the isotropic cycles associated to these points are mapped to each other by the orthogonal reflection $\mathfrak{R}$ defined by the non-isotropic cycle $p$.

Recall that an isotropic cycle associated to the point $\infty$ is the 0-cycle $C$ (where $C$ is any non-zero number in the coordinate field), and an isotropic cycle associated to a finite point ***m*** is any scalar multiple of the 2-cycle $\boldsymbol{X}.\boldsymbol{X} - 2\boldsymbol{m}.\boldsymbol{X} + \boldsymbol{m}.\boldsymbol{m}$.

The following proposition provides another characterization of such conjugate points.

**Proposition 6.1:** *Two distinct points* ***m*** *and* ***m'*** *are conjugate with respect to p if and only if the non-isotropic cycle p belongs to the linear subspace spanned by the isotropic cycles associated with* ***m*** *and* ***m'***.

**Proof.** Let $q$ and $q^*$ be two isotropic cycles associated with ***m*** *and* ***m'*** respectively. If ***m*** and ***m'*** are conjugate with respect to $p$, we must have:

$$\lambda q^* = \mathfrak{R}(q) = q - 2\,\frac{\langle q,p \rangle}{\langle p,p \rangle}\, p \quad \text{for some non-zero scalar } \lambda$$

This means that $q^*$ is generated by $q$ and $p$, i.e., the projective classes of $p$, $q$ and $q^*$ are collinear. The point ***m'*** (represented by $q^*$) is in the line through ***m*** (represented by $q$) and the projective class of $p$.

Because ***m*** and ***m'*** are distinct, the coefficient of $p$ in the above equation must also be non-zero. Accordingly, to say $q^*$ is generated by $q$ and $p$ is the same as saying $p$ is generated by $q$ and $q^*$.

Conversely, suppose that a non-isotropic cycle $p$ is generated by $q$ and $q^*$. The pencil $\hat{L}$ generated by $p$ and $q$ intersects $V$ in exactly two distinct points, ***m*** and ***m'***. The conjugate of ***m*** with respect to the cycle $p$ belongs to the same pencil, as proven above. Such a conjugate point of course lies in $V$, and hence that conjugate point must be the same as point $\boldsymbol{m'}$ because $V$ is stable under the orthogonal reflection defined by $p$. ■

**<u>Proposition 6.2:</u>** *Let u and v be two non-isotropic cycles in F. Recall our convention that the zero set of a 1-cycle includes also the point at infinity.*

*(a) If the zero set Z(v) passes through two points* **m** *and* **m'** *in* $E \cup \{\infty\}$ *that are conjugate with respect to u, then* $<u, v> = 0$.

*(b) Conversely, assume that* $<u, v> = 0$. *Then for any point* **m** *in the zero set Z(v), its conjugate point* **m'** *with respect to u also belongs to Z(v).*

**Proof**. Let $q$ and $q^*$ be two isotropic cycles associated with the points $\boldsymbol{m}$ and $\boldsymbol{m'}$ in the extended affine space $E \cup \{\infty\}$. If Z($v$) passes through $\boldsymbol{m}$ and $\boldsymbol{m'}$, that means $v$ is orthogonal to both $q$ and $q^*$ in light of Proposition 2.1. Therefore $v$ is orthogonal to all cycles spanned by $q$ and $q^*$, which includes $u$ by virtue of Proposition 6.1.

Conversely, if $v$ is orthogonal to $u$, and $v$ is also orthogonal to the isotropic 2-cycle $q$, then $v$ must be orthogonal to all cycles spanned by $u$ and $q$, which includes $q^*$. That means $v(\boldsymbol{m'}) = 0$ or that the conjugate point $\boldsymbol{m'}$ also lies on the zero set $Z(v)$ just like $\boldsymbol{m}$ does. ■

We now discuss briefly the well-known stereographic projection mapping and its generalization to our case. Let $L$ be the vector space of elements of the form $(x, \boldsymbol{y}, z)$, where $x$ and $z$ are numbers in $K$, and $\boldsymbol{y}$ is a vector in the space $E$.

Define a scalar product $< _ , _>$ on $L$ as follows. For $t = (x, \boldsymbol{y}, z)$ and $t^* = (x^*, \boldsymbol{y}^*, z^*)$, we define the scalar product $<t, t^*>$ of $t$ and $t^*$ as the number:

$$<t, t^*> = \boldsymbol{y}.\boldsymbol{y}^* + xx^* - zz^*$$

This scalar product on $L$ is clearly non-degenerate, since it makes $L$ isometric to the sum of $E$ and an Artinian plane. We will call this scalar product the Lorentz product by virtue of its analogy with the classical Lorentz metric.

As with our foregoing study of $F$ and its scalar product, we can look at the projective space $\hat{L}$ and the set $U$ of all isotropic elements in the space $\hat{L}$. [5]

Now consider the following linear mapping $\mathcal{S}$ from $L$ to $F$:

$$\mathcal{S}\colon\ (\mathrm{x}, \mathbf{y}, \mathrm{z}) \ \rightarrow\ \frac{(\mathrm{z}-\mathrm{x})}{2}\,\mathbf{X.X} \ -\ \mathbf{y.X}\ +\frac{(\mathrm{z}+\mathrm{x})}{2}$$

$\mathcal{S}$ is clearly an invertible linear mapping of $L$ onto $F$. When $x \neq z$, it maps $(x, \boldsymbol{y}, z)$ to a 2-cycle centered at $\boldsymbol{y}/(z-x)$. It maps $(1, \mathbf{0}, 1)$ to the constant function $p(\boldsymbol{X}) = 1$.

Moreover, $\mathcal{S}$ is compatible with the inner products of the two spaces in the sense that:

$$<u, v> \text{ (Lorentz product on } L) \ = \ <\mathcal{S}(u), \mathcal{S}(v)> \text{ (cycle product on } F)$$

To see this, note that by virtue of bilinear conditions, it is sufficient to verify

$$<u, u> \ = \ <\mathcal{S}(u), \mathcal{S}(u)> \quad \text{for any } u$$

If $u = (x, \boldsymbol{y}, z)$, we have $<u, u> \ = \ \boldsymbol{y.y} + x^2 - z^2$.

On the other side, for $\mathcal{S}(u) = \dfrac{(z-x)}{2}\,\boldsymbol{X.X} \ -\ \boldsymbol{y.X}\ + \dfrac{(z+x)}{2}$, we have:

$$<\mathcal{S}(u), \mathcal{S}(u)> = \ \boldsymbol{y.y} - 4\,\frac{(z-x)}{2}\,\frac{(z+x)}{2} = \boldsymbol{y.y} - (z^2 - x^2) = \boldsymbol{y.y} + x^2 - z^2$$

---

[5] In the case of a finite-dimensional space $E$ over an ordered field $K$, the set $U$ lies entirely in the affine subset of the projective space $\hat{L}$ corresponding to $\mathrm{z} = 1$ because the equation $x^2 + \boldsymbol{y.y} = 0$ has no non-zero solution. If we identify the vector $(x, \boldsymbol{y}, 1)$ in $L$ with the point $(x, \boldsymbol{y})$ in $K^{1+n}$, then the set $U$ can be identified with the unit sphere $S^n$ of all points $(x, \boldsymbol{y})$ in $k^{1+n}$ satisfying the equation $x^2 + \boldsymbol{y.y} = 1$.

Consequently, $\mathcal{S}$ is an isometry of the vector space $L$ (endowed with the Lorentz product) onto the vector space $F$ (endowed with the cycle product). The isometry $\mathcal{S}$ induces a natural bijection between the set $U$ of isotropic elements in the projective space $\hat{L}$ and the set $V$ of isotropic elements in the projective space $\hat{F}$, under which the point $(1, \mathbf{0}, 1)$ of $U$ is mapped to the point $\infty$ of $V = E \cup \{\infty\}$ and each other point $(x, \mathbf{y}, z)$ $(x \neq z)$ of $U$ is mapped to the point $\mathbf{y}/(z - x)$ in $V = E \cup \{\infty\}$ .

This is exactly the good old stereographic projection when $E$ is a Euclidean space. From our current point of view, this mapping represents a trace of a larger isometry $\mathcal{S}$ between the spaces $L$ and $F$. From this perspective, intersections of $U$ with hyperplanes in the space $\hat{L}$ are transformed into intersections of $V = E \cup \{\infty\}$ with hyperplanes in the space $\hat{F}$ in a linear and orthogonal fashion.

## **Section 7. Inversive geometry in dimension one**

The results of this paper apply to any field $K$ of characteristic $\neq 2$ in case of dimension 1, because multiplication in the field $K$ is naturally an anisotropic bilinear form.

In this one-dimensional case, the vector space $F$ of cycles has linear dimension $1 + 2 = 3$ and the projective space $\hat{F}$ has projective dimension $1 + 1 = 2$ (a projective plane). The elements of $F$ are expressions $Ax^2 + Bx + C$ in one variable $x$, where $A$, $B$, $C$ vary in $K$.

The set $V$ of isotropic cycles is identified with the extended line $K \cup \{\infty\}$. The constant expressions $C$ for 0-cycles represent the exceptional point $\infty$, while the expressions $A(x - u)^2$ represent the finite point $u$ for any number $u$ in the field $K$, i.e., the isotropic 2-cycles centered at the point $u$.

Any 1-cycle $Bx + C$ has exactly one zero in $K$, namely the number $(- C/B)$. Geometrically, the projective hyperplane orthogonal to that 1-cycle intersects $V$ in exactly two points, namely at $\infty$ and at $(- C/B)$. The affine reflection defined by 1-cycle $Bx + C$ is the transformation that sends $\infty$ to itself, and any number $x$ to $x' = x - 2\,\dfrac{Bx+C}{B} = -x - 2(C/B)$. Such a transformation is a combination of the symmetry $x \to -x$ and the translation $x \to x - 2(C/B)$.

For a 2-cycle $p = Ax^2 + Bx + C$, rewrite it as $A(x + B/2A)^2 - As$ where $s = (B^2 - 4AC)/4A^2$. The 2-cycle $p$ is non-isotropic when the equation $Ax^2 + Bx + C = 0$ does not have a double root. That means each non-isotropic 2-cycle either has two distinct roots or none at all. Geometrically, the orthogonal complement of $\hat{p}$ either intersects the set V in two distinct points, or not at all.

The inversion relative to such a non-zero 2-cycle is the transformation $x \to x^*$ described by the equation: $(x^* + B/2A)(x + B/2A) = (B^2 - 4AC)/4A^2$. If $B = 0$ and $A = 1$, i.e, for a non-isotropic 2-cycle centered at 0, that inversion is simply the mapping $x^* = -(C/x)$.

From the basic background outlined above, it follows that the inversive transformations in this case are simply all the fractional linear transformations described by:

$$x \to x' = \frac{ax+b}{cx+d} \quad \text{where } a, b, c, d \text{ are numbers such that } ad - bc \neq 0,$$

because such transformations are generated by the translations and inversions described above.

Therefore, the group of all fractional linear transformations of $V = K \cup \{\infty\}$ can be described as the projective transformation group $PGL(1, K)$ of the projective line $\mathbf{P}(K^2)$ if we identify $K \cup \{\infty\}$ with $\mathbf{P}(K^2)$ under the standard identification that maps $K$ to the affine set $(x, 1)$ and $\infty$ to the projective point corresponding to the line $(x, 0)$ in $K^2$.

Accordingly, in the one-dimensional case, inversive transformations can either be identified with fractional linear transformations on a projective line, or with orthogonal projective transformations in the space $F$ of cycles. This double connection with projective transformations in dimension 1 (which is valid for any field $K$ of characteristic $\neq 2$) can give us some new insights on the classical geometry of the projective line and interesting generalizations. Here we can do no better than refer the readers to [2] and [3] for more information.

## Section 8. A Natural Cross Ratio

Consider any non-degenerate quadratic space $F$ of Witt index 1 (over a field $K$ of characteristic $\neq 2$). The Witt index 1 here means that any isotropic subspace of $F$ has dimension at most 1.

The geometric space of interest to us is the projective quadric $Q$ of $F$ under the action of the projective orthogonal group $G$. Following Felix Klein, we are interested in finding numbers associated with configurations of points in the projective quadric $Q$ that are invariant under the action of the projective orthogonal group $G$.

If $u$ and $v$ are any representative vectors of two distinct points in the projective quadric $Q$, the pairing $\langle u, v\rangle$ is invariant under the action of $G$. However, that pairing by itself cannot be an invariant, because different representative vectors will give us different values of the pairing.

Note that $\langle u, v\rangle$ is always non-zero, because otherwise the two vectors $u$ and $v$ will generate a totally isotropic subspace of dimension 2, which by hypothesis is not possible. With that in mind, we can consider the following cross ratio associated to any four distinct points in the quadric $Q$, with representative vectors $x, y, z, t$:

$$\langle x, z\rangle.\langle y, t\rangle / \langle x, t\rangle.\langle y, z\rangle$$

This cross ratio is well-defined because any pairing value of two distinct points is always non-zero as noted above, and moreover the cross ratio gives us the same value regardless of which representative vectors are chosen for the four distinct points due to cross cancellation. This cross ratio is naturally invariant under the action of $G$. The same invariant applies if we just work with a subset of the quadric $Q$ which is stable under the action of a subgroup of $PO(F)$.

In particular, consider the quadratic space $F$ associated with each anisotropic quadratic space $E$ as constructed in this paper. The quadratic space $F$ has Witt index 1, and therefore we have a natural cross ratio for the corresponding projective quadric $V$ under the action of the projective orthogonal group $PO(F)$.

In the case where the anisotropic space $E$ is the one-dimensional vector space $K$ with multiplication (see Section 7), the projective quadric $V$ can be naturally identified with the projective line $K \cup \{\infty\}$ and the projective orthogonal group $PO(F)$ can be identified with the projective transformation group $PGL(1, K)$. For two finite points $A$ and $B$ on the projective line, if we take the 2-cycles $x^2 - 2Ax + A^2$ and $x^2 - 2Bx + B^2$ as vector representatives of these two points, then their cycle pairing is $-2(A - B)^2$. That means the natural cross ratio in this case is equal to the square of the classical cross ratio of a projective line.

That fact suggests the following proposition.

**Proposition 8.1:** *Any transformation of a projective line over a field K of characteristic ≠ 2 that leaves invariant the square of the classical cross ratio must be a projective transformation.*

**Proof**. Let $p$ be a transformation of the projective line $K \cup \{\infty\}$ that leaves invariant the square of the classical cross ratio. We will show that if 0, 1 and ∞ are fixed by $p$, then $p$ must be the identity transformation. This will prove the proposition, because any three distinct points can be mapped by a projective transformation to 0, 1 and ∞.

For any point $\lambda \neq 0$, 1 and ∞, the classical cross ratio $[\lambda, 1; 0, \infty] = \lambda$. So $p(\lambda) = \lambda$ or $-\lambda$. This implies that $p(-1) = -1$, as we already have $p(1) = 1$.

Let $\lambda \neq 0, 1, -1$ and ∞. Consider the cross ratio $[\lambda, -1; 1, \infty] = (1 - \lambda)/2$. If $p(\lambda) = -\lambda$, then the square of that cross ratio must be equal to the square of the cross ratio $[-\lambda, -1; 1, \infty] = (1 + \lambda)/2$. That means $(1 - \lambda)^2 = (1 + \lambda)^2$ or $-2\lambda = 2\lambda$, a contradiction as $\lambda \neq 0$. Accordingly, $p(\lambda) = \lambda$ for all $\lambda$, i.e., $p$ is the identity transformation. ■

It happens that some well-known geometric spaces, such as the hyperbolic plane or hyperbolic spaces of higher dimensions, can be regarded as subset of a quadric $Q$ in a space $F$ of Witt index 1, and the relevant geometric transformation group can be identified with a subgroup of $PO(F)$. In those situations, the cross ratio described here gives a natural invariant for the geometry. In addition, if we restrict attention to certain subgroups of $PO(F)$, other invariants may also come up in a natural way. For example, if we look at the geometry defined by the subgroup of orthogonal transformations leaving invariant a given nonisotropic vector $p$ of $F$, then the expression

$$\langle x, p\rangle.\langle y, p\rangle / \langle x, y\rangle$$

clearly gives us another natural invariant for that geometry. I hope that the interested reader will have a happy time exploring these questions.